\documentclass[12pt]{article}

\usepackage{subcaption}

\usepackage{float}

\usepackage{multirow}

\usepackage{authorindex} 

\usepackage{tikz}

\usetikzlibrary{calc,backgrounds,arrows,matrix}

\usepackage{enumerate}

\usepackage{color}
\definecolor{lightblue}{rgb}{0,0.2,0.5}

\usepackage[colorlinks=true, urlcolor=blue,linkcolor=blue, citecolor=lightblue]{hyperref}

\usepackage[round,sort,comma,numbers]{natbib}

\bibpunct{\textcolor{lightblue}{(}}{\textcolor{lightblue}{)}}{,}{a}{}{;}

\usepackage{amssymb,amsmath}

\usepackage{graphicx}

\usepackage{graphics}
\usepackage{color} 

\DeclareMathAlphabet{\eufrak}{U}{}{}{}
\SetMathAlphabet\eufrak{normal}{U}{euf}{m}{n}
\SetMathAlphabet\eufrak{bold}{U}{euf}{b}{n}

\allowdisplaybreaks

 \def\qu{{\mathord{\mathbb Z}}}

 \def\Var{{\mathrm{{\rm Var}}}}

 \def\sZZ{{\rm Z\kern-.45em{}Z}}

 \def\sQQ{{\kern 0.27em \vrule height1.45ex width0.03em depth0em
           \kern-0.30em \rm Q}}
 \def\qu{{\mathchoice
         {\sQQ}
         {\sQQ}
   {\kern 0.225em \vrule height1.05ex width0.025em depth0em \kern-0.25em \rm Q}
   {\kern 0.180em \vrule height0.78ex width0.020em depth0em \kern-0.20em \rm Q}
         }}
 \def\sGG{{\kern 0.27em \vrule height1.45ex width0.03em depth0em
           \kern-0.30em \rm G}}
 \def\gg{{\mathchoice
         {\sGG}
         {\sGG}
   {\kern 0.225em \vrule height1.05ex width0.025em depth0em \kern-0.25em \rm G}
   {\kern 0.180em \vrule height0.78ex width0.020em depth0em \kern-0.20em \rm G}
         }}

 \newtheorem{prop}{Proposition}[section]

 \newtheorem{corollary}[prop]{Corollary}
 
 \newtheorem{remark}[prop]{Remark}

\numberwithin{equation}{section}

 \def\P{{\mathord{\mathbb P}}}

\newcommand{\re}{\mathrm{e}}

 \newcounter{hyp}
\usepackage{aeguill}

\newenvironment{Proof}{\removelastskip\par\medskip \noindent{\em Proof.} \rm}{\penalty-20\null\hfill$\square$\par\medbreak}

\def\bprf{\begin{Proof}}
\def\nprf{\end{Proof}}
\def\bdes{\begin{description}}
\def\ndes{\end{description}}

\newtheorem{thm}{Theorem}[section]

\def\bdef{\begin{defn}}
\def\ndef{\end{defn}}
\def\bthm{\begin{thm}}
\def\nthm{\end{thm}}
\def\bprop{\begin{prop}}
\def\nprop{\end{prop}}
\def\brmk{\begin{remark}}
\def\nrmk{\end{remark}}
\def\bexa{\begin{exa}}
\def\nexa{\end{exa}}
\def\blem{\begin{lem}}
\def\nlem{\end{lem}}
\def\bcor{\begin{cor}}
\def\ncor{\end{cor}}
\def\bexe{\begin{exe}}
\def\nexe{\end{exe}}

\newcommand{\E}{\mathbb{E}}

\newcommand{\real}{\mathbb{R}}

\def\Var{\mathop{\hbox{\rm Var}}\nolimits}

\def\og{\leavevmode\raise.3ex
     \hbox{$\scriptscriptstyle\langle\!\langle$~}}
\def\fg{\leavevmode\raise.3ex
     \hbox{~$\!\scriptscriptstyle\,\rangle\!\rangle$}~}

 \title{\Huge
   Recursive computation of the Hawkes cumulants
 }

 \author{Nicolas Privault
\\ 
\small
Division of Mathematical Sciences 
\\ 
\small 
School of Physical and Mathematical Sciences 
\\ 
\small
Nanyang Technological University 
\\ 
\small 
21 Nanyang Link 
\\ 
\small
Singapore 637371
}

\begin{document}

\maketitle 

\vspace{-0.9cm}

\baselineskip0.6cm
 
\begin{abstract} 
 We propose a recursive method for the computation of the cumulants
 of self-exciting point processes of Hawkes type, 
 based on standard combinatorial tools such as Bell polynomials. 
 This closed-form approach is easier to implement 
 on higher-order cumulants 
 in comparison with existing methods based on
 differential equations, tree enumeration or
 martingale arguments.
 The results are corroborated by Monte Carlo simulations,
 and also apply to the computation of joint cumulants
 generated by multidimensional self-exciting processes. 
\end{abstract} 
 
 
\baselineskip0.7cm

\section{Introduction}
Hawkes processes were introduced in
\cite{hawkes1971} as self-exciting point processes
representing an alternative to doubly stochastic point processes. 
In recent years they have found applications in many fields,
from neuroscience, see e.g. \cite{cardanobile},
to genomics analysis, see e.g. \cite{reynaud-bouret},
as well as finance \cite{embrechts2} and social media \cite{younglee}. 
As noted in \cite{jovanovic}, the analysis
of statistical properties of Hawkes processes is still incomplete,
in particular in terms of moments, cumulants and other
statistical parameters
 such as skweness and kurtosis. 

\medskip 

In \cite{dassios-zhao2} the moment and probability generating functions of
(generalized) Hawkes processes and their intensity have been obtained by ODE methods, with the computation
of first and second moments in the stationary case, see also
 \cite{errais}. 
 In \cite{bacry}, a stochastic calculus and martingale approach has been
 applied to the computation of first and second moments, however it
 seems difficult to generalize to higher orders.
 In \cite{jovanovic}, a tree-based method for the computation of cumulants
 has been introduced, with an explicit computation of third order cumulants.
 However, this type of algorithm requires to perform tree enumerations. 

\medskip

Third-order cumulant expressions for Hawkes processes
 have been used in \cite{achab} for the
estimation of branching ratio matrices in the analysis of
order books,
and in \cite{ocker},
\cite{montangie} 
for the estimation of third order correlations in
spiking neuronal networks.
Higher order cumulants can also be
useful in order to provide finer estimates of the
evolution of time correlations and
 of the probability density functions of neuronal membrane potentials
 by Gram-Charlier density expansions, see e.g.
 \cite{brigham-destexhe}, \cite{neuron}. 

\medskip 

In this paper, we derive a general recursion formula
using the standard Bell polynomials for the computation of the cumulants of  
 a self-exciting point process on
$\real^d$, $d\geq 1$, with immigrant intensity $\nu (dx)$ and
branching intensity $\gamma (dx)$ on $(\real^d,{\cal B}(\real^d ))$.
Our approach is based on a recursive relation for
the Probability Generating Functional (PGFl) $G_z$ of a self-exciting
point process from a single point at $z\in \real^d$,
derived in Proposition~\ref{fjklds}.
Such an implicit relation has already been observed in 
\cite{adamopoulos}, and applied in e.g. 
\cite{bordenave-torrisi} to large deviations,
however 
it does not seem to have been exploited for the computation of cumulants. 

\medskip

 In Section~\ref{s1} we start by reviewing
 the combinatorial approach of \S~3.2 of \cite{consul} to the computation
 of the cumulants of the integer-valued Borel distribution,
 and show that it can be extended as an
 explicit recursion using Bell polynomials.
 This provides an elementary model for subsequent computations, as
 the Borel distribution can be used to represent the cardinality of
 a self-exciting Poisson cluster point process. 

 \medskip

 Next, in Section~\ref{s2} we extend this argument to the computation of
 the cumulants of self-exciting Hawkes Poisson cluster processes
 in Proposition~\ref{fjklds}, with an extension to
 the computation of joint cumulants. 
 This provides a closed-form alternative, suitable
 for systematic higher-order computations,
 to the tree-based approach of \cite{jovanovic}. 
 Explicit computations for the time-dependent 
 third and fourth cumulants and skewness and kurtosis
 of Hawkes processes with exponential kernels are presented
 in Section~\ref{s5},
 and are confirmed by Monte Carlo estimates. 

\subsubsection*{Cumulants, Fa\`a di Bruno formula and Bell polynomials}
We close this section with background results on combinatorics
 that will be needed in the sequel. 
Recall that if $f(t)$ admits the formal series expansion
$$
f(t) = \sum_{n=1}^\infty \frac{a_n}{n!}t^n,
$$ 
by the Fa\`a di Bruno formula we have
\begin{eqnarray}
 \label{dfjkl0} 
\re^{f ( t )} - 1 
 & = & 
 \sum_{n=1}^\infty
 \frac{t^n}{n!}
 \sum_{k=1}^n
 \frac{n!}{k!}
 \sum_{l_1+\cdots + l_k = n\atop l_1,\ldots , l_k \geq 1}
 \frac{a_{l_1} \cdots a_{l_k}}{l_1!\cdots l_k!}
 \\
 \label{dfjkl}
  & = &
 \sum_{n=1}^\infty
 \frac{t^n}{n!}
 B_n ( a_{l_1}, \ldots ,a_{l_k} ) 
\end{eqnarray} 
 where the sum \eqref{dfjkl0} holds on the integer compositions
 $(l_1,\ldots ,l_k)$ of $n$,
 see e.g. Relation~(2.5) in \cite{elukacs},
$$ 
 B_n ( a_1 , \ldots , a_n ) 
  = 
 \sum_{k=1}^n
 B_{n,k} ( a_1 , \ldots , a_{n-k+1} )
=
 \sum_{k=1}^n
 \frac{n!}{k!} 
 \sum_{l_1+\cdots + l_k=n \atop 
 l_1\geq 1,\ldots ,l_k \geq 1 
 } 
 \frac{a_{l_1}}{l_1!} 
 \cdots 
 \frac{a_{l_k}}{l_k!} 
$$ 
 is the complete Bell polynomial of degree $n \geq 1$, 
 and
 $$ 
 B_{n,k} ( a_1 , \ldots , a_{n-k+1} ) 
=
 \frac{n!}{k!} 
 \sum_{l_1+\cdots + l_k=n \atop 
 l_1\geq 1,\ldots ,l_k \geq 1 
 } 
 \frac{a_{l_1}}{l_1!} 
 \cdots 
 \frac{a_{l_k}}{l_k!},
 \qquad 1\leq k \leq n, 
$$ 
 is the partial Bell polynomial 
 of order $(n,k)$. 
 Given the Moment Generating Function (MGF) 
\aimention{Thiele, T.N.}
\begin{equation} 
\nonumber 
 M_X (t) : = 
 \E \big[ \re^{tX} \big] 
 = 
 1 + \sum_{n\geq 1} 
 \frac{t^n}{n!} 
 \E [ X^n ], 
 \qquad 
 t \in \real, 
\end{equation} 
of a random variable $X$,
 the {cumulants} of $X$ are the coefficients 
 $\big(\kappa_X^{(n)}\big)_{n\geq 1}$ appearing in the log-MGF series expansion 
\begin{equation} 
\label{cgf} 
\log M_X (t)
 =
\log \big( \E\big[\re^{tX}\big] \big) 
 =
 \sum_{n\geq 1} \kappa_X^{(n)} \frac{t^n}{n!}, 
 \quad 
 t \in \real. 
\end{equation} 
 The moments $\E[X^n]$ of a random variable $X$ 
 are linked to its cumulants $\big(\kappa_X^{(n)}\big)_{n\geq 1}$
 through the relation 
\begin{equation} 
\nonumber 
 \E [ X^n ] 
 = 
 \sum_{k=1}^n 
 \sum_{\pi_1 \cup \cdots \cup \pi_k = \{ 1, \ldots , n \}} 
 \kappa_X^{(|\pi_1|)} \cdots \kappa_X^{(|\pi_k|)}, 
\end{equation} 
 which runs over the partitions 
 $\pi_1,\ldots , \pi_a$ of the set $\{ 1 , \ldots , n \}$,
 where $|\pi_i|$ denotes the cardinality of $\pi_i$. 
 By the Fa\`a di Bruno formula, \eqref{cgf} can be inverted as
\begin{equation} 
\nonumber 
 \kappa_X^{(n)} = \sum_{k=1}^n 
 (k-1)! (-1)^{k-1} 
 \sum_{\pi_1 \cup \cdots \cup \pi_k = \{ 1,\ldots ,n \} } 
 \E \big[ X^{|\pi_1|} \big] \cdots \E \big[ X^{|\pi_k|} \big],
 \quad n \geq 1, 
\end{equation} 
 see e.g. Theorem~1 of \cite{elukacs}, 
 and also \cite{leonov}, 
 Relations~(2.8)-(2.9) in \cite{mccullagh}, or 
 Corollary~5.1.6 in \cite{stanley}. 
 The third and fourth cumulants can be used to define the skewness
 $\kappa_X^{(3)}/\big( \kappa_X^{(2)} \big)^{3/2}$
 and the excess kurtosis $\kappa_X^{(4)}/ \big( \kappa_X^{(2)}\big)^2$
 of $X$. 
 \section{Borel cumulants} 
\label{s1}
In this section we consider the recursive computation
of the cumulants of
integer-valued Borel-distributed random variables
using the Fa\`a di Bruno formula.
For this, we review the method of
\S~3.2 of \cite{consul} which applies to 
  {Lagrangian distributions},
  and note that it admits an explicit formulation
  using Bell polynomials.    
 Let $(X_n)_{n\geq 0}$ be a branching process started at $X_0=1$
 with Poisson distributed offspring count $N$ of parameter $\mu \in (0,1)$. 
Denoting by $X$ the total count of offsprings generated by $(X_n)_{n\geq 0}$
and letting $\big(X^{(l)}\big)_{l\geq 1}$ denote a sequence of independent copies of $X$,
the Probability Generating Function (PGF) of $X$ can be estimated 
by the standard branching recursion 
\begin{eqnarray} 
  \nonumber
  G_X (s) & = & \E[ s^X ] 
\\ 
\nonumber
   & = & 
 s \E \left[ \prod_{l=1}^N s^{X^{(l)}} \right]
 \\
 \nonumber
  & = & 
 s \sum_{k\geq 0} 
 \E \left[ \prod_{l=1}^k s^{X^{(l)}}
 \right] 
 \P ( N = k ) 
\\ 
\nonumber
   & = & 
 s \re^{-\mu} \sum_{n\geq 0} 
 \big( \E \big[ s^{X^{(1)}} \big] \big)^n 
 \frac{\mu^n}{n!} 
\\ 
\nonumber 
& = & 
 s G_\mu ( G_X ( s ) ), 
 \qquad 
 -1 \leq s \leq 1, 
\end{eqnarray} 
 where 
$$
 G_\mu (s) := \re^{-\mu} \sum_{n=0}^\infty \frac{\mu^n}{n!} s^n
 = \re^{ \mu ( s-1)}, \qquad s\in [-1,1], 
 $$
 is the PGF of the Poisson distribution with mean $\mu>0$.
 The equation
\aimention{Haight, F.A.}
\aimention{Breuer, M.A.}
\begin{equation}
\label{pgf} 
 G_X (s) = s G_\mu ( G_X ( s ) ), 
 \qquad -1 \leq s \leq 1, 
\end{equation} 
 can be solved using Lagrange series,  
 see page~145 of \cite{polya-szego}, 
 showing that $X$ has the {Borel distribution} 
\index{Lagrangian distribution}
\index{distribution!Lagrangian}
\index{Borel distribution}
\index{distribution!Borel}
\aimention{Consul, P.C.}
\aimention{Famoye, F.}
\aimention{P{\'o}lya, G.} 
\aimention{Szeg{\"o}, G.} 
 $$
 \P ( X = n )
 = \re^{-\mu n}\frac{(\mu n)^{n-1}}{n!},
 \qquad n \geq 1,
 $$
 which belongs to the class
 of {Lagrangian distributions}, see \S~8.4 of \cite{consul}. 
 The following proposition then extends the relations
 (3.12) in \cite{consul} 
 for the computation of the cumulants of the Borel
 distribution, via a general expression based on the Bell polynomials. 
 Another, less direct, recursion can be found in 
 \S~8.4.3 in \cite{consul}, based on the derivatives
 of moments of $X$ with respect to $\mu$.
 \begin{prop} Let $X$ be a Borel distributed random variable with parameter
   $\mu \in (0,1)$. 
 We have $\kappa^{(1)}_X = 1/(1-\mu)$ and the induction relation 
$$ 
 \kappa_X^{(n)}
 =
 \frac{\mu}{1-\mu}
 \big( B_n \big( \kappa_X^{(1)}, \ldots , \kappa_X^{(n)} \big)
 - \kappa_X^{(n)} \big)
=
 \frac{\mu}{1-\mu}
  \sum_{k=2}^n
 B_{n,k} \big( \kappa_X^{(1)}, \ldots , \kappa_X^{(n-k+1)} \big) 
, 
 \quad n\geq 2, 
 $$ 
 where $B_n$, resp. $B_{n,k}$,
 is the complete, resp. partial, Bell polynomial. 
\end{prop}
\begin{Proof} 
From \eqref{pgf} the moment generating function $M_X (t) = \E [ \re^{tX}]
= G(\re^t)$ satisfies
$$
\log M_X (t) = t + \mu ( M_X (t) -1) = t + \mu ( \re^{\log M_X (t)} -1), \quad
t\in \real, 
$$
 see also Relation~(19) in \cite{haight} with $r=1$. 
 Based on the cumulant expansion \eqref{cgf} 
 and the Fa\`a di Bruno formula \eqref{dfjkl}, we have 
\begin{eqnarray*} 
 \sum_{n\geq 1} \kappa_X^{(n)} \frac{t^n}{n!}
 & = & 
 \log M_X  (t) 
 \\
  & = & t + \mu \big( \re^{\log M_X  (t )} - 1 \big) 
 \\
 & = & ( \mu + \kappa_X^{(1)} ) t +
 \mu
 \sum_{n=2}^\infty
 \frac{t^n}{n!}
 B_n \big( \kappa_X^{(1)}, \ldots , \kappa_X^{(n)} \big), 
\end{eqnarray*} 
 which shows that 
 $
 \kappa_X^{(1)} 
  = 
 1 
 +
 \mu \kappa_X^{(1)}
 $
  and 
  $$
  \kappa_X^{(n)} 
  =
  \mu
  B_n \big( \kappa_X^{(1)}, \ldots , \kappa_X^{(n)} \big)
  = 
  \mu \kappa_X^{(n)} 
  +
  \sum_{k=2}^n
 B_{n,k} \big( \kappa_X^{(1)}, \ldots , \kappa_X^{(n-k+1)} \big),
  \qquad n \geq 2. 
$$
\end{Proof}
\noindent
In particular, since $B_2(x_1,x_2)=x_1^2+x_2$ we have 
$$ 
 \kappa_X^{(2)}
 = 
 \frac{\mu}{1-\mu}
 \left( B_2 \left( \frac{1}{1-\mu}, \kappa_X^{(2)} \right)
 - \kappa_X^{(2)} \right)
 = 
 \frac{\mu}{(1-\mu)^3}. 
$$ 
 Given that $B_3(x_1,x_2,x_3)=x_1^3+3x_1x_2+x_3$, 
 we have 
$$
 \kappa_X^{(3)}
  = 
 \frac{\mu}{1-\mu}
 \left( B_3 \left( \frac{1}{1-\mu},
 \frac{\mu}{(1-\mu)^3} ,
 \kappa_X^{(3)} \right)
 - \kappa_X^{(3)} \right)
 =
 \mu \frac{1+2\mu}{(1-\mu)^5}
. 
$$ 
 Since 
 $B_4(x_1,x_2,x_3,x_4)=x_1^4+6x_1^2x_2+4x_1x_3+3x_2^2+x_4$,
 we find 
$$ 
 \kappa_X^{(4)}
 = 
 \frac{\mu}{1-\mu}
 \left( B_4 \left( \frac{1}{1-\mu},\frac{\mu}{(1-\mu)^3},
 \mu \frac{1+2\mu}{(1-\mu)^5}, \kappa_X^{(4)} \right)
 - \kappa_X^{(4)} \right)
 = \mu \frac{ 1 + 8\mu + 6 \mu^2 }{(1-\mu)^7}, 
$$ 
 which recovers (8.85) page 159 of \cite{consul}. 
\section{Hawkes cumulants} 
\label{s2}
In this section we work in the cluster process framework of \cite{hawkes}.
We consider a self-exciting point process on $\real^d$, $d\geq 1$,
with Poisson
offspring intensity $\gamma (dx)$ and Poisson
immigrant intensity $\nu (dx)$ on $\real^d$, built on the space
$$
 \Omega = \big\{
 \xi = \{ x_i \}_{i\in I} \subset \real^d \ : \
 \#( A \cap \xi ) < \infty 
 \mbox{ for all compact } A\in {\cal B} (\real^d ) 
 \big\}
 $$
 of locally finite configurations on $\real^d$, whose elements 
 $\xi \in \Omega$ are identified with the Radon point measures 
 $\displaystyle \xi (dz) = \sum_{x\in \xi} \epsilon_x (dz)$, 
 where $\epsilon_x$ denotes the Dirac measure at $x\in \real^d$. 
 In particular, any initial immigrant point $z \in \real^d$ branches into a Poisson
 random sample $\xi_\gamma (\cdot + dz)$ centered at $z$, with intensity
 measure $\gamma (\cdot + dz)$. 
 We let
$$ 
 G_z (f) = f(z) \E\left[ \prod_{x\in \xi} f(z+x) \right] 
$$ 
 denote the Probability Generating Functional (PGFl) of the branching process starting 
 from a single point at $z\in \real^d$,
 for sufficiently integrable $f:\real^d\to \real$. 
The next proposition states a recursive property
for the Probability Generating Functional $G_z(f)$, see
also Theorem~1 in \cite{adamopoulos}. 
\begin{prop}
  \label{fjklds}
  The Probability Generating Functional $G_z(f)$ satisfies
$$ 
 G_z (f) 
 = f(z) \exp \left( \int_{\real^d} ( G_{z+x} ( f ) - 1 ) \gamma ( dx) \right),
 \qquad z\in \real^d, 
$$ 
 and the PGFl of the Hawkes process with immigrant intensity $\nu (dz)$ is given by
$$ 
  G_\nu (f)
  = \exp \left( \int_{\real^d} ( G_z (f) - 1 ) \nu ( dz ) \right). 
$$
\end{prop}
\begin{Proof} 
  Viewing the self-exciting point process $\xi$ as a marked
  point process we have, see e.g. Lemma~6.4.VI of \cite{daley}, 
\begin{eqnarray*} 
 G_z (f) 
 & = & f(z) \E\left[ \prod_{x\in \xi} f(z+x) \right] 
\\
 & = & f(z) \E\left[ \prod_{x\in \xi_\gamma} \left( \prod_{y\in \xi} f(z+x+y) \right) \right] 
\\
 & = & f(z) \E\left[ \prod_{x\in \xi_\gamma} \E\left[ \prod_{y\in \xi} f(z+x+y) \right] \right] 
\\
 & = & f(z) \E\left[ \prod_{x\in \xi_\gamma} G_{z+x} ( f ) \right] 
\\
& = & \re^{-\gamma (\real^d) } f(z) \sum_{k=0}^\infty \frac{1}{k!}
\int_{(\real^d)^k} G_{z+x_1}(f) \cdots G_{z+x_k} (f) \gamma (dx_1)\cdots \gamma ( dx_k) 
\\
& = & f(z) \exp \left( \int_{\real^d} ( G_{z+x} ( f ) - 1 ) \gamma ( dx) \right), 
\end{eqnarray*} 
 and 
\begin{eqnarray*} 
  G_\nu (f)
  & = &
  \re^{-\nu (\real^d )} \sum_{n=0}^\infty \frac{1}{n!}
  \int_{(\real^d )^n} G_{z_1}(f) \cdots G_{z_n}(f) \nu ( dz_1 ) \cdots \nu ( dz_n ) 
  \\
    & = &
  \exp \left( \int_{\real^d} ( G_z (f) - 1 ) \nu ( dz ) \right). 
\end{eqnarray*} 
\end{Proof}
\noindent
 Let
$$ 
 M_z (f) = G_z \big(e^f\big) = \E\left[ \exp \left( f(z) +
   \sum_{x\in \xi} f(z+x) \right) \right] 
$$ 
 denote the Moment Generating Functional (MGFl) of the
 stochastic integral
 $\displaystyle \sum_{x\in \xi} f(x)$ 
 given that the cluster process $\xi$ starts 
 from a single point at $z\in \real^d$. 
 The following corollary is an immediate consequence of
 Proposition~\ref{fjklds}, see also Proposition~2.6 in
 \cite{bogachev2} for Poisson cluster processes.
\begin{corollary} 
  The Moment Generating Functional $M_z(f)$ satisfies
  the recursive relation 
 \begin{equation}
   \label{mfdsf0} 
 M_z (f) 
 = \exp \left( f(z) + \int_{\real^d} ( M_{z+x} ( f ) - 1 ) \gamma ( dx) \right),
 \qquad z\in \real^d.
\end{equation} 
 The MGFl of the Hawkes process with immigrant intensity $\nu (dz)$ is given by
 \begin{equation}
   \label{mfdsf} 
 M_\nu (f)
 = 
 \exp \left( \int_{\real^d} ( M_z (f) - 1 ) \nu ( dz ) \right). 
\end{equation} 
\end{corollary}
\noindent 
The next proposition provides a way to compute the
cumulants $\kappa_z^{(n)}(f)$ of $\displaystyle \sum_{x\in \xi} f(x)$
 by an induction relation based on the Bell polynomials. 
 Note that the sum of coefficients in $B_n (x_1,\ldots , x_n)$ is the Bell 
 number 
$$
 B_n = \sum_{k=1}^n
 \frac{1}{k!}
 \sum_{l_1+\cdots + l_k = n\atop l_1,\ldots , l_k \geq 1}
 \frac{n!}{l_1!\cdots l_k!}
 $$
 that represents the count of partitions of a set of $n$ elements.  
 In the sequel we consider the integral operator $\Gamma$ defined as
 $$
 \Gamma f (z) = \int_{\real^d} f(z+y ) \gamma (dy), \qquad z\in \real^d, 
 $$
 and the inverse operator $(I_d-\Gamma)^{-1}$ given by 
 $$
 (I_d-\Gamma)^{-1} f(z ) 
   = f(z) + \sum_{m=1}^\infty \int_{\real^d} \cdots \int_{\real^d} f(z+x_1+\cdots + x_m) \gamma ( dx_1)
  \cdots \gamma ( dx_m), \quad z\in \real^d. 
  $$
   
\begin{prop}
\label{p1}
 The first cumulant
 $\kappa_z^{(1)}(f)$ of $\displaystyle \sum_{x\in \xi} f(x)$
 given that $\xi$ is started 
 from a single point at $z\in \real^d$ 
 is given by 
 $\kappa_z^{(1)}(f)  = (I_d-\Gamma )^{-1} f(z)$ for $n=1$, 
 and for $n\geq 2$ by the induction relation 
\begin{eqnarray*} 
 \kappa_z^{(n)}(f) 
 & = & 
 ( I_d - \Gamma )^{-1} \Gamma 
 \big( B_n \big( \kappa_{z+\cdot }^{(1)}, \ldots , \kappa_{z+ \cdot }^{(n)} \big)
 - \kappa_{z+\cdot}^{(n)} \big)
 \\
  & = & 
  \sum_{k=2}^n
  ( I_d - \Gamma )^{-1} \Gamma 
  B_{n,k}
   \big( \kappa_{z+\cdot }^{(1)}, \ldots , \kappa_{z+\cdot}^{(n-k+1)} \big). 
\end{eqnarray*} 
\end{prop}
\begin{Proof} 
  By \eqref{cgf}, \eqref{mfdsf0}
 and the Fa\`a di Bruno formula \eqref{dfjkl}, we have
 \begin{eqnarray}
   \nonumber 
  \sum_{n=1}^\infty \frac{t^n}{n!} \kappa_z^{(n)}(f)
& = &
  \log M_z ( t f) 
  \\
   \nonumber 
   & = & tf(z) + \int_{\real^d} \big( \re^{\log M_{z+x} ( tf )} - 1 \big) \gamma ( dx)
\\
\label{fjlkdsf} 
 & = & tf(z)
 +
 t 
 \int_{\real^d}
 \kappa_{z+x}^{(1)}
 \gamma ( dx)
 +
 \sum_{n=2}^\infty
 \frac{t^n}{n!}
 \int_{\real^d}
 B_n \big( \kappa_{z+x}^{(1)}, \ldots , \kappa_{z+x}^{(n)} \big)
  \gamma ( dx),  
\end{eqnarray} 
 hence
 \begin{eqnarray*} 
  \kappa_z^{(1)}(f) & = & f(z) + \int_{\real^d} \kappa_{z+x}^{(1)}(f) \gamma ( dx)
  \\
  & = & f(z) + \sum_{m=1}^\infty \int_{\real^d} \cdots \int_{\real^d} f(z+x_1+\cdots + x_m) \gamma ( dx_1) \cdots \gamma ( dx_m), 
\end{eqnarray*} 
 as solution of the renewal equation
$$
 \kappa_z^{(1)}(f)  = h(z) + \int_{\real^d} \kappa_{z+x}^{(1)}(f)  \gamma (dx),
 \quad z \in \real^d. 
$$
 For $n\geq 2$, \eqref{fjlkdsf} yields 
$$ 
 \kappa_z^{(n)}(f) 
  = 
 \int_{\real^d}
  B_n \big( \kappa_{z+x}^{(1)}, \ldots , \kappa_{z+x}^{(n)} \big)
  \gamma ( dx)
 = 
 \Gamma \kappa_{z+\cdot }^{(n)}(f) 
  +
 \Gamma 
 \big( B_n \big( \kappa_{z+\cdot }^{(1)}, \ldots , \kappa_{z+\cdot }^{(n)} \big)
 - \kappa_{z+\cdot }^{(n)} \big),  
$$ 
 or
 $$ 
 (I_d - \Gamma ) \kappa_z^{(n)}(f) 
  = 
 \Gamma 
 \big( B_n \big( \kappa_{z+\cdot }^{(1)}, \ldots , \kappa_{z+\cdot }^{(n)} \big)
 - \kappa_{z+\cdot }^{(n)} \big), 
$$ 
 which yields 
\begin{eqnarray*} 
  \lefteqn{
    \kappa_z^{(n)}(f) 
 = 
 ( I_d - \Gamma )^{-1} \Gamma 
 \big( B_n \big( \kappa_{z+\cdot }^{(1)}, \ldots , \kappa_{z+\cdot }^{(n)} \big)
 - \kappa_{z+\cdot }^{(n)} \big)
  }
  \\
 & = & 
 \sum_{m=1}^\infty
  \int_{\real^d} \cdots \int_{\real^d} \big( B_n \big( \kappa_{z+x_1+\cdots + x_m}^{(1)}, \ldots , \kappa_{z+x_1+\cdots + x_m}^{(n)} \big)
 - \kappa_{z+x_1+\cdots + x_m}^{(n)} \big)
 \gamma ( dx_1) \cdots \gamma ( dx_m ) 
 \\
 & = & 
 \sum_{m=1}^\infty
  \sum_{k=2}^n
  \int_{\real^d} \cdots \int_{\real^d} 
  B_{n,k} \big( \kappa_{z+x_1+\cdots + x_m}^{(1)}, \ldots , \kappa_{z+x_1+\cdots + x_m}^{(n-k+1)} \big)
 \gamma ( dx_1) \cdots \gamma ( dx_m ) 
 , \qquad n \geq 2. 
\end{eqnarray*} 
\end{Proof}
Unconditional cumulants can be obtained in the next corollary as a consequence
of Proposition~\ref{p1}. 
\begin{corollary} 
\label{c1} 
 The cumulant of order $n\geq 2$ of $\displaystyle \sum_{x\in \xi} f(x)$
 is given by the recursion 
$$
    \kappa^{(n)} (f)
  = 
\int_{\real^d}
B_n \big( \kappa_z^{(1)}(f) , \ldots , \kappa_z^{(n)}(f) \big) 
\nu ( dz),
$$
 with 
\begin{eqnarray*} 
    B_n \big( \kappa_z^{(1)}(f) , \ldots , \kappa_z^{(n)}(f) \big)
  & = & 
    ( I_d - \Gamma )^{-1}
    \big( B_n \big( \kappa_{z+\cdot }^{(1)}, \ldots , \kappa_{z+\cdot }^{(n)} \big)
    - \kappa_{z+\cdot }^{(n)} \big)
    \\
     & = & 
       \sum_{k=2}^n
       ( I_d - \Gamma )^{-1}
       B_{n,k} \big( \kappa_{z+\cdot }^{(1)}, \ldots , \kappa_{z+\cdot }^{(n-k+1)} \big),
       \quad z\in \real^d.
\end{eqnarray*} 
\end{corollary} 
\begin{Proof} 
\noindent
 By \eqref{cgf}, \eqref{mfdsf}
 and the Fa\`a di Bruno formula \eqref{dfjkl}, we have
\begin{eqnarray*} 
 \sum_{n=1}^\infty \frac{t^n}{n!} \kappa^{(n)}(f)
& = &
 \log M_\nu (f)
 \\
 & = &  
 \int_{\real^d} ( M_z (f) - 1 ) \nu ( dz )
 \\
  & = &  
 \int_{\real^d} ( \re^{\log M_z (f)} - 1 ) \nu ( dz )
  \\
  & = &  
 \sum_{n=1}^\infty
 \frac{t^n}{n!}
 B_n \big( \kappa_z^{(1)}(f) , \ldots , \kappa_z^{(n)}(f) \big)
 \nu ( dz), 
\end{eqnarray*}
and therefore
$$ 
\kappa^{(n)} (f) =
\int_{\real^d}
B_n \big( \kappa_z^{(1)}(f) , \ldots , \kappa_z^{(n)}(f) \big) 
 \nu ( dz), \qquad n \geq 2.
 $$ 
 We conclude from the equalities
\begin{eqnarray*}
  \lefteqn{
    B_n \big( \kappa_z^{(1)}(f) , \ldots , \kappa_z^{(n)}(f) \big)
  }
  \\
  & = &
  \kappa_z^{(n)}(f)
  + \big( 
  B_n \big( \kappa_z^{(1)}(f) , \ldots , \kappa_z^{(n)}(f) \big)
  -
  \kappa_z^{(n)}(f) \big) \big)
  \\
  & = &
 \sum_{m=1}^\infty
  \int_{\real^d} \cdots \int_{\real^d} \big( B_n \big( \kappa_{z+x_1+\cdots + x_m}^{(1)}, \ldots , \kappa_{z+x_1+\cdots + x_m}^{(n)} \big)
 - \kappa_{z+x_1+\cdots + x_m}^{(n)} \big)
 \gamma ( dx_1) \cdots \gamma ( dx_m ) 
 \\
 & &
 +
  \big( 
  B_n \big( \kappa_z^{(1)}(f) , \ldots , \kappa_z^{(n)}(f) \big)
  -
  \kappa_z^{(n)}(f) \big) \big)
  \\
  & = &
 \sum_{m=0}^\infty
  \int_{\real^d} \cdots \int_{\real^d} \big( B_n \big( \kappa_{z+x_1+\cdots + x_m}^{(1)}, \ldots , \kappa_{z+x_1+\cdots + x_m}^{(n)} \big)
 - \kappa_{z+x_1+\cdots + x_m}^{(n)} \big)
 \gamma ( dx_1) \cdots \gamma ( dx_m ), 
\end{eqnarray*}
 that follow from Proposition~\ref{p1}.    
\end{Proof} 
\subsubsection*{Second cumulant}
\noindent
 For $n=2$, Proposition~\ref{p1} shows that 
$$ 
 \kappa_z^{(2)} (f) 
 =
 \sum_{m=1}^\infty 
 \int_{\real^d} \cdots \int_{\real^d} \big( \kappa^{(1)}_{z+x_1+\cdots + x_m}(f)\big)^2 \gamma ( dx_1) \cdots \gamma ( dx_m), 
$$
 and by Corollary~\ref{c1} we have 
\begin{eqnarray*} 
  \kappa^{(2)} (f)
  & = & 
\int_{\real^d} \kappa_z^{(2)}(f) \nu ( dz)
  +
 \int_{\real^d} \big( \kappa_z^{(1)}(f) \big)^2 \nu ( dz)
 \\
 & = & 
  \sum_{m=0}^\infty 
 \int_{\real^d} \cdots \int_{\real^d} \big( \kappa^{(1)}_{z+x_1+\cdots + x_m}(f)\big)^2 \gamma ( dx_1) \cdots \gamma ( dx_m) \nu ( dz), 
\end{eqnarray*} 
see e.g. Proposition~2 in \cite{bacry} and
 Eq.~(37) in \cite{jovanovic}.

\subsubsection*{Third cumulant}
 For $n=3$, we have
\begin{eqnarray} 
  \nonumber
  \kappa_z^{(3)} (f) 
& = &  
 3 \sum_{m=1}^\infty 
 \int_{\real^d} \cdots \int_{\real^d} \kappa_{z+x_1+\cdots + x_m}^{(1)} (f) \kappa_{z+x_1 + \cdots + x_m }^{(2)} (f)
 \gamma ( dx_1) \cdots \gamma ( dx_m ) 
 \\
 \label{ef3}
 & & + 
 \sum_{m=1}^\infty
 \int_{\real^d} \cdots \int_{\real^d} \big( \kappa_{z+x_1+\cdots + x_m}^{(1)} (f) \big)^3 
 \gamma ( dx_1) \cdots \gamma ( dx_m ) 
,
\end{eqnarray} 
 and 
\begin{eqnarray} 
  \label{k3}
  \lefteqn{
  \! \! \! \!    \kappa^{(3)} (f) 
 = 
 \int_{\real^d} B_3 \big(
 \kappa_z^{(1)},\kappa_z^{(2)},\kappa_z^{(3)}\big) \nu (dz) 
  }
  \\
 \nonumber
 & = &  
 \int_{\real^d} \big( \kappa_z^{(1)} \big)^3 \nu (dz) 
 +
 3 \int_{\real^d} \kappa_z^{(1)}\kappa_z^{(2)} \nu (dz) 
 +
 \int_{\real^d} \kappa_z^{(3)} \nu (dz) 
 \\
\label{f1} 
 & = &  
 3 \sum_{m=0}^\infty 
 \int_{\real^d} \cdots \int_{\real^d} \kappa_{z+x_1+\cdots + x_m}^{(1)} (f) \kappa_{z+x_1 + \cdots + x_m }^{(2)} (f)
 \gamma ( dx_1) \cdots \gamma ( dx_m ) \nu (dz)  
 \\
\label{f2} 
 & & + 
 \sum_{m=0}^\infty
 \int_{\real^d} \cdots \int_{\real^d} \big( \kappa_{z+x_1+\cdots + x_m}^{(1)} (f) \big)^3 
 \gamma ( dx_1) \cdots \gamma ( dx_m )  \nu (dz), 
\end{eqnarray} 
which corresponds to Eq.~(39) in \cite{jovanovic}.
\subsubsection*{Fourth cumulant}
 For $n=4$, we have
\begin{eqnarray*} 
    \kappa_z^{(4)} (f) 
 & = & 
 6 
 \sum_{m=1}^\infty
 \int_{\real^d} \cdots \int_{\real^d}
 \big( \kappa_{z+x_1+\cdots + x_m}^{(1)} (f) \big)^2 
 \kappa_{z+x_1 + \cdots + x_m }^{(2)} (f)
 \gamma ( dx_1) \cdots \gamma ( dx_m ) 
 \\
 & & + 
4 
 \sum_{m=1}^\infty
 \int_{\real^d} \cdots \int_{\real^d} \kappa_{z+x_1+\cdots + x_m}^{(1)} (f) \kappa_{z+x_1 + \cdots + x_m }^{(3)} (f) \gamma ( dx_1) \cdots \gamma ( dx_m ) 
     \\
 & & + 
 3
 \sum_{m=1}^\infty
 \int_{\real^d} \cdots \int_{\real^d} \big( \kappa_{z+x_1+\cdots + x_m}^{(2)} (f) \big)^2
 \gamma ( dx_1) \cdots \gamma ( dx_m ) 
 \\
 & & + 
 \sum_{m=1}^\infty
 \int_{\real^d} \cdots \int_{\real^d}
 \big( \kappa_{z+x_1+\cdots + x_m}^{(1)} (f) \big)^4
 \gamma ( dx_1) \cdots \gamma ( dx_m ) 
 ,
\end{eqnarray*} 
 and 
\begin{eqnarray} 
  \label{k4}
  \kappa^{(4)} (f) 
& = &  
 \int_{\real^d} B_4 \big(
 \kappa_z^{(1)},\kappa_z^{(2)},\kappa_z^{(3)},\kappa_z^{(4)}\big) \nu (dz) 
 \\
\nonumber 
 & = &  
 \int_{\real^d} \big( \kappa_z^{(1)} \big)^4 \nu (dz) 
 +
 6 \int_{\real^d} \big( \kappa_z^{(1)} \big)^2\kappa_z^{(2)} \nu (dz) 
 \\
\nonumber 
  & & +
 4 \int_{\real^d} \kappa_z^{(1)}\kappa_z^{(3)} \nu (dz) 
 +
 3 \int_{\real^d} \big( \kappa_z^{(2)} \big)^2 \nu (dz) 
 +
 \int_{\real^d} \kappa_z^{(4)} \nu (dz) 
 \\
\label{ff1} 
  & = & 
 6 
 \sum_{m=0}^\infty
 \int_{\real^d} \cdots \int_{\real^d}
 \big( \kappa_{z+x_1+\cdots + x_m}^{(1)} (f) \big)^2 
 \kappa_{z+x_1 + \cdots + x_m }^{(2)} (f)
 \gamma ( dx_1) \cdots \gamma ( dx_m ) 
 \\
\nonumber 
 & & + 
4 
 \sum_{m=0}^\infty
 \int_{\real^d} \cdots \int_{\real^d} \kappa_{z+x_1+\cdots + x_m}^{(1)} (f) \kappa_{z+x_1 + \cdots + x_m }^{(3)} (f) \gamma ( dx_1) \cdots \gamma ( dx_m ) 
     \\
\nonumber 
 & & + 
 3
 \sum_{m=0}^\infty
 \int_{\real^d} \cdots \int_{\real^d} \big( \kappa_{z+x_1+\cdots + x_m}^{(2)} (f) \big)^2
 \gamma ( dx_1) \cdots \gamma ( dx_m ) 
 \\
\label{ff4} 
 & & + 
 \sum_{m=0}^\infty
 \int_{\real^d} \cdots \int_{\real^d}
 \big( \kappa_{z+x_1+\cdots + x_m}^{(1)} (f) \big)^4
 \gamma ( dx_1) \cdots \gamma ( dx_m ) 
 .
\end{eqnarray} 

\noindent
We note that the count of $4$ terms in \eqref{f1}-\eqref{f2} 
and the total count of $6+4\times 4+3\times 1 + 1\times 1=24$
 terms in \eqref{ff1}-\eqref{ff4} match the ones obtained in Figure~4 of
\cite{jovanovic} using tree enumeration. 
\subsubsection*{Joint cumulants} 
 The expression of Proposition~\ref{p1} can be
extended to joint cumulants by standard combinatorial
arguments. 
\begin{prop}
\label{djklds}
 For $n\geq 2$, the joint cumulants 
 $\kappa_z^{(n)}(f_1,\ldots , f_n)$
 of $\displaystyle \sum_{x\in \xi} f_1(x), \ldots , \sum_{x\in \xi} f_n(x)$
 given that $\xi$ is started 
 from a single point at $z\in \real^d$ 
 are given by the induction relation 
$$ 
 \kappa_z^{(n)}(f_1,\ldots , f_n) 
 = 
 \sum_{k=2}^n
 \sum_{\pi_1\cup \cdots \cup \pi_k = \{1,\ldots , n\}} 
 \sum_{m=1}^\infty
 \int_{\real^d} \cdots \int_{\real^d}
 \prod_{j=1}^k
 \kappa_{z+x_1+\cdots + x_m}^{(|\pi_j|)} ((f_i)_{i\in \pi_j})
 \gamma ( dx_1) \cdots \gamma ( dx_m ) 
,
$$ 
$n \geq 2$, where the above sum is over set partitions
$\pi_1\cup \cdots \cup \pi_k=\{1,\ldots , n\}$, $k=2,\ldots , n$.
\end{prop}
 As in Corollary~\ref{c1}, we obtain the expressions 
\begin{eqnarray*} 
  \lefteqn{
    \kappa^{(n)} (f_1,\ldots , f_n) = 
  \sum_{k=1}^n
 \sum_{\pi_1\cup \cdots \cup \pi_k = \{1,\ldots , n\}} 
 \int_{\real^d} 
 \prod_{j=1}^k
 \kappa_z^{(|\pi_j|)}((f_i)_{i\in \pi_j})
 \nu ( dz ) 
  }
 \\
  & = & 
 \sum_{k=2}^n
 \sum_{\pi_1\cup \cdots \cup \pi_k = \{1,\ldots , n\}} 
 \sum_{m=0}^\infty
 \int_{\real^d} \cdots \int_{\real^d}
 \prod_{j=1}^k \kappa_{z+x_1+\cdots + x_m}^{(|\pi_j|)}((f_i)_{i\in \pi_j})
 \gamma ( dx_1) \cdots \gamma ( dx_m ) \nu (dz),
\end{eqnarray*} 
 as a consequence of Proposition~\ref{djklds}. 
 \subsubsection*{Second joint cumulant}
 \noindent
  We have 
$$ 
 \kappa_z^{(2)} (f_1,f_2) 
 =
 \sum_{m=1}^\infty 
 \int_{\real^d} \cdots \int_{\real^d}
 \kappa^{(1)}_{z+x_1+\cdots + x_m} (f_1)
 \kappa^{(1)}_{z+x_1+\cdots + x_m} (f_2) 
 \gamma ( dx_1) \cdots \gamma ( dx_m), 
$$
 and 
 $$ 
\kappa^{(2)} (f_1,f_2) =
\int_{\real^d} \kappa_z^{(2)}(f_1,f_2) \nu ( dz)
  +
 \int_{\real^d} \kappa_z^{(1)}(f_1) \kappa_z^{(1)}(f_2)
 \nu ( dz).
$$ 
\subsubsection*{Third joint cumulant}
 For $n=3$, we have
\begin{eqnarray} 
  \nonumber
\lefteqn{ 
  \! \!
  \kappa_z^{(3)} (f_1,f_2,f_3) 
 = 
 \sum_{m=1}^\infty 
 \int_{\real^d} \cdots \int_{\real^d} \kappa_{z+x_1+\cdots + x_m}^{(1)} (f_1)
 \kappa_{z+x_1 + \cdots + x_m}^{(2)} (f_2,f_3)
 \gamma ( dx_1) \cdots \gamma ( dx_m ) 
}
\\
\nonumber
 &  & 
 + \sum_{m=1}^\infty 
 \int_{\real^d} \cdots \int_{\real^d} \kappa_{z+x_1+\cdots + x_m}^{(1)} (f_2)
 \kappa_{z+x_1 + \cdots + x_m}^{(2)} (f_1,f_3)
 \gamma ( dx_1) \cdots \gamma ( dx_m ) 
 \\
\nonumber
  & & 
 + \sum_{m=1}^\infty 
 \int_{\real^d} \cdots \int_{\real^d} \kappa_{z+x_1+\cdots + x_m}^{(1)} (f_3)
 \kappa_{z+x_1 + \cdots + x_m}^{(2)} (f_1,f_2)
 \gamma ( dx_1) \cdots \gamma ( dx_m ) 
 \\
\nonumber 
 & & 
 + 
 \sum_{m=1}^\infty
 \int_{\real^d} \cdots \int_{\real^d}
 \kappa_{z+x_1+\cdots + x_m}^{(1)} (f_1)
 \kappa_{z+x_1+\cdots + x_m}^{(1)} (f_2)
 \kappa_{z+x_1+\cdots + x_m}^{(1)} (f_3)
 \gamma ( dx_1) \cdots \gamma ( dx_m ),
\end{eqnarray} 
 and
\begin{eqnarray*} 
  \lefteqn{
    \kappa^{(3)} (f_1,f_2,f_3) 
 = 
 \int_{\real^d}
 \kappa_z^{(1)}(f_1)
 \kappa_z^{(1)}(f_2)
 \kappa_z^{(1)}(f_3)
 \nu (dz) 
  }
  \\
 & &
 +
 \int_{\real^d} \kappa_z^{(1)}(f_1)\kappa_z^{(2)}(f_2,f_3) \nu (dz) 
 +
 \int_{\real^d} \kappa_z^{(1)}(f_2)\kappa_z^{(2)}(f_1,f_3) \nu (dz) 
 +
 \int_{\real^d} \kappa_z^{(1)}(f_3)\kappa_z^{(2)}(f_1,f_2) \nu (dz) 
 \\
 & &
 +
 \int_{\real^d} \kappa_z^{(3)}(f_1,f_2,f_3) \nu (dz).
\end{eqnarray*} 
 Similar expressions for $\kappa^{(4)} (f)$ can be obtained from \eqref{k4}. 
\section{Example - exponential kernel} 
\label{s5} \noindent 
In this section we take $d=1$ and consider the exponential kernel
$\gamma (dx) = a {\bf 1}_{[0,\infty )} (x) \re^{-bx} dx$, $0< a < b$, 
and constant Poisson intensity $\nu (dz) = \nu dz$, $\nu >0$. In this case,
$\displaystyle N_t (\xi ) := \xi ( [0,t]) = \sum_{x\in \xi} {\bf 1}_{[0,t]}(x)$ 
 defines the self-exciting Hawkes process with stochastic intensity 
$$ 
 \lambda_t =
 a \int_0^t \re^{-b(t-s)} dN_s, \qquad t\in \real_+.
$$
 The recursive calculation of the cumulants
 $\kappa^{(n)}_z (t):=\kappa^{(n)}_z ({\bf 1}_{[0,t]})$
 will be performed using 
 the family of functions $e_{p,\eta}(x) := x^p \re^{\eta x} {\bf 1}_{[0,t]}(x)$,
 $\eta < b$, $p\geq 0$, which satisfy the relation 
\begin{eqnarray*} 
 ( I_d - \Gamma )^{-1} \Gamma 
 e_{p,\eta} ( z ) 
 & = & 
\sum_{n=1}^\infty
 \int_0^t \cdots \int_0^t
 e_{p,\eta} ( z+x_1+\cdots + x_n )
 \gamma (dx_1) \cdots \gamma (dx_n )
 \\
  & = &  
a 
 \int_0^{t-z}
 y^p \re^{(\eta + a -b)y}
 dy, \qquad z\in [0,t], 
\end{eqnarray*} 
 with
 \begin{eqnarray} 
\nonumber 
     ( I_d - \Gamma )^{-1} \Gamma 
     e_{0,\eta} ( z) 
     & = & 
     a \re^{\eta z}
     {\bf 1}_{(-\infty , t]}(z)
          \frac{\re^{(\eta + a -b) (t-z)} - 1}{\eta + a-b}
 \\
   \label{c1jk} 
  & = & a \re^{\eta z} \frac{\re^{(\eta + a -b) t}e_{0,-\eta + b- a }(z) - e_{0,0}(z)}{\eta + a-b}, 
\end{eqnarray} 
 and 
 \begin{equation}
     \label{c1jk2} 
       ( I_d - \Gamma )^{-1} \Gamma 
     e_{1,\eta} ( z) 
   = 
 a \re^{\eta z} \frac{e_{0,0}(z)
   + \re^{ (\eta + a - b )t  }
   e_{0,-\eta + b - a }(z) 
   ( (\eta + a - b )(t-z) -1)}{(\eta + a - b )^2}, 
\end{equation} 
 where $ ( I_d - \Gamma )^{-1} \Gamma
     e_{p,\eta} ( z) 
$ can be similarly evaluated for $p \geq 2$. 
\subsubsection*{First cumulant} 
 We have 
\begin{eqnarray}
\nonumber 
 \kappa^{(1)}_z (t)  
 & = &
 1 + \sum_{n=1}^\infty
 \int_0^t \cdots \int_0^t
 e_{0,0} ( z+x_1+\cdots + x_n )
 \gamma ( dx_1) \cdots \gamma ( dx_n ) 
 \\
 \label{e1}
 & = &
 1 + \frac{a}{a-b} \re^{(a-b )t}e_{0,b-a} (z) - \frac{a}{a-b} e_{0,0}(z) 
, \qquad z\in \real_+, 
\end{eqnarray} 
 which recovers 
$$ 
 \E [ N_t] = \kappa^{(1)} (t) 
 = 
 \int_0^t \kappa_z^{(1)} (t) \nu ( dz ) 
 =
 \frac{\nu}{(b-a)^2}
 \big( - a
 + b(b-a) t 
   +a \re^{ (a-b) t} 
 \big)
, 
$$ 
 as solution of the differential equation
$$
d\E[N_t] = \nu + a \int_0^t \re^{-b  ( t-s)} d\E [N_s].
$$ 
\subsubsection*{Second cumulant} 
 Using \eqref{c1jk}, we have
\begin{eqnarray} 
  \nonumber
  \lefteqn{
 \kappa^{(2)}_z (t)  = 
 \sum_{m=1}^\infty 
 \int_0^t \cdots \int_0^t
 \big( \kappa^{(1)}_{z+x_1+\cdots + x_m} (t) \big)^2 
 \gamma ( dx_1) \cdots \gamma ( dx_m)
  }
  \\
\nonumber
  & = &
 \frac{b^2}{(b-a)^2}
 \sum_{m=1}^\infty
 \int_0^t \cdots \int_0^t
 e_{0,0}(z ) 
 \gamma ( dx_1 ) \cdots \gamma ( dx_m ) 
  \\
\nonumber
  & &
- \frac{2a b}{(b-a)^2}
\re^{(b-a)t}
\sum_{m=1}^\infty
 \int_0^t \cdots \int_0^t
 e_{0,b-a}(z+x_1+\cdots + x_m ) 
 \gamma ( dx_1 ) \cdots \gamma ( dx_m ) 
 \\
\nonumber
   & & + 
\frac{a^2}{(b-a)^2}
\re^{2(b-a)t}
\sum_{m=1}^\infty
 \int_0^t \cdots \int_0^t
 e_{0,2(b-a)}(z+x_1+\cdots + x_m )  
 \gamma ( dx_1 ) \cdots \gamma ( dx_m ) 
\\
\nonumber 
& = &
\frac{ a b ^2
}{(b-a)^3}
 ( 1 - \re^{ - (b-a)  ( t-z)} )
- \frac{2a^2 b  }{(b-a)^2}
(t-z) \re^{ - (b-a) (t-z)}
- a^3 \frac{\re^{-2 (b-a)  (t-z)} - \re^{ - (b-a) (t-z)}}{(b-a)^3}, 
\\
\label{e2}
\end{eqnarray} 
and 
\begin{eqnarray*} 
\int_0^t \big( \kappa_z^{(1)}(t) \big)^2 \nu (dz) 
& = & 
\int_0^t  \left( 
  \frac{b }{b-a}
  -      \frac{a}{b-a}
\re^{ - (b-a) (t-z)}
\right)^2  
\nu (dz)
\\
 & = & 
 \frac{\nu b ^2 t }{(b-a)^2}
- \frac{2 \nu b  a}{(b-a)^3}
( 1 - \re^{ - (b-a) t})
+ \frac{\nu a^2}{2(b-a)^3}
( 1 - \re^{-2 (b-a) t}),
\end{eqnarray*} 
 hence 
\begin{align*} 
 &     \Var [ N_t] = \kappa^{(2)} (t) 
  \\
  & = 
 -\frac{\nu}{2 (a - b)^4}
 \Big(
 6 a b^2 -a^2 b + 2 b^3 ( a - b ) t
 + 2 a \big( a^2 
 - 3 b^2 + 2 ab( a - b ) t
 \big) \re^{(a - b) t}
 + a^2 ( b - 2 a ) \re^{2 (a - b) t}
 \Big). 
\end{align*} 

\noindent
 The following figures are plotted with $\nu=1$, $a=0.5$, $b=1$, and $10^7$
 Monte Carlo samples.
 
\begin{figure}[H]
  \centering
 \begin{subfigure}[b]{0.49\textwidth}
    \includegraphics[width=1\linewidth, height=5cm]{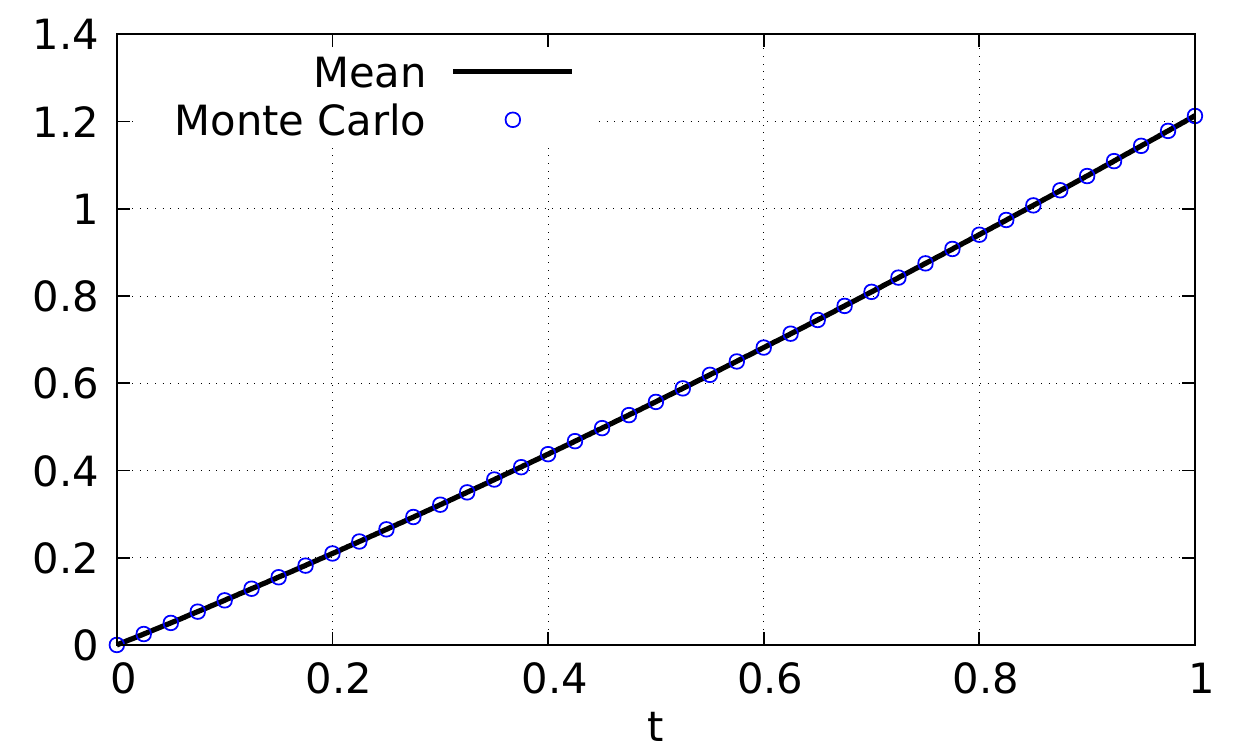}
    \caption{First cumulant $\kappa^{(1)} (t)$.} 
 \end{subfigure}
  \begin{subfigure}[b]{0.49\textwidth}
    \includegraphics[width=1\linewidth, height=5cm]{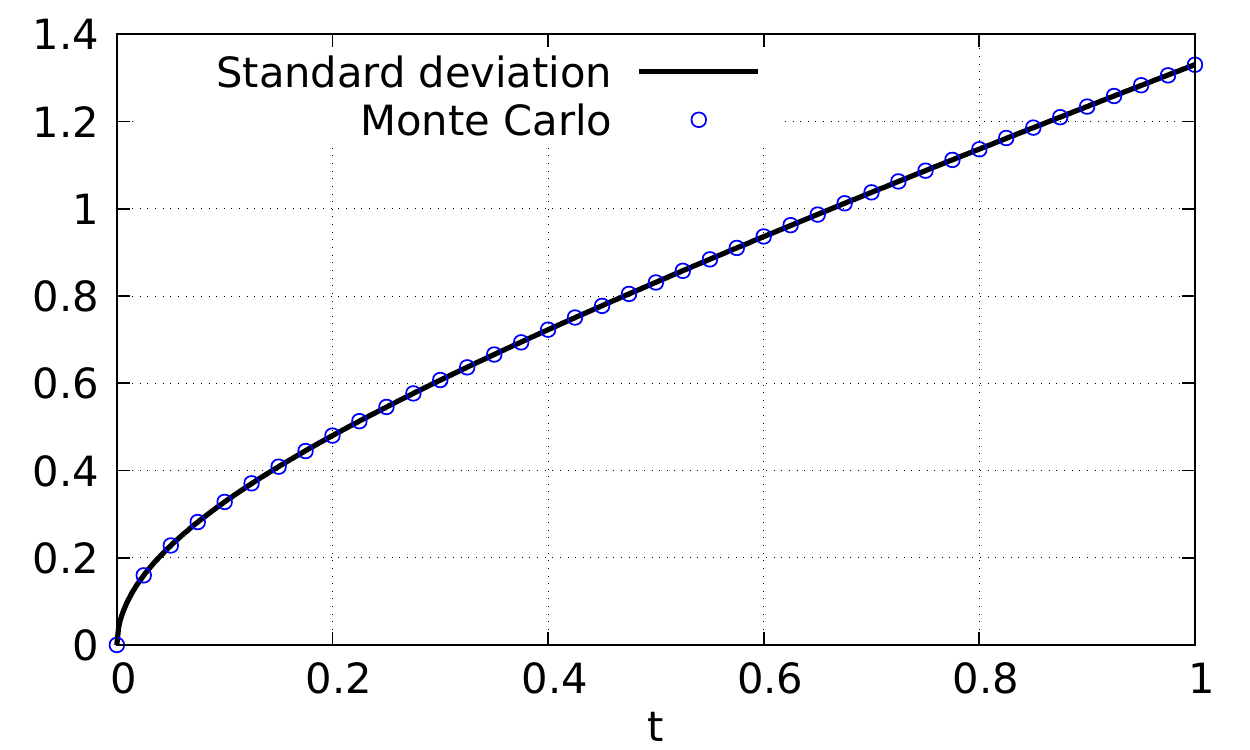} 
    \caption{Square root of second cumulant $\kappa^{(2)} (t)$.} 
  \end{subfigure}
    \caption{Mean and standard deviation with exponential kernel.} 
\end{figure}
\subsubsection*{Third cumulant} 
The recursive computation of $\kappa^{(3)} (t)$ can be carried out
from \eqref{ef3}-\eqref{k3} and \eqref{c1jk}-\eqref{c1jk2} using Mathematica 
 based on 
$\kappa_{z+x_1+\cdots + x_m}^{(1)} (f)$, 
$\kappa_{z+x_1+\cdots + x_m}^{(2)} (f)$
given in \eqref{e1}-\eqref{e2},
which yields 
\begin{align*} 
  & \kappa^{(3)}(t) = 
  -\frac{\nu }{6 (a-b)^6}
  \Big(
     42 a b^4 
    +30 a^2 b^3 
    -7 a^3 b^2 
    +a^4 b
    +
    6 b^4 (
    2 a^2 - a b 
    - b^2 
    ) t 
    \\
    &
         +
    3 \big( 
    18 a^3 b^2 
    -16 a^2 b^3
    - a^4 b 
    -14 a b^4
    -2 a^5
    +
    6 a^2 b ( 
    4 a b^2 
    -4 b^3
    + a^2 b 
    - a^3 ) t
 - 6 a^3 b^2 ( a-b )^2 t^2 
    \big) \re^{ (a- b) t}
    \\
    &
    +
    9 \big(
    2 a^2 b^3 
    -5 a^3 b^2 
    - a^4 b 
    +2 a^5
    +
    2 a^3 b ( b^2 -3 a b  
    +2 a^2 ) t 
    \big) \re^{ 2(a-b) t }
        -
    a^3 \big( 2 b^2 
    -11 a b 
    +12 a^2 \big)
    \re^{3 (a-b) t}
    \Big).
\end{align*}

\begin{figure}[H]
  \centering
 \begin{subfigure}[b]{0.49\textwidth}
    \includegraphics[width=1\linewidth, height=5cm]{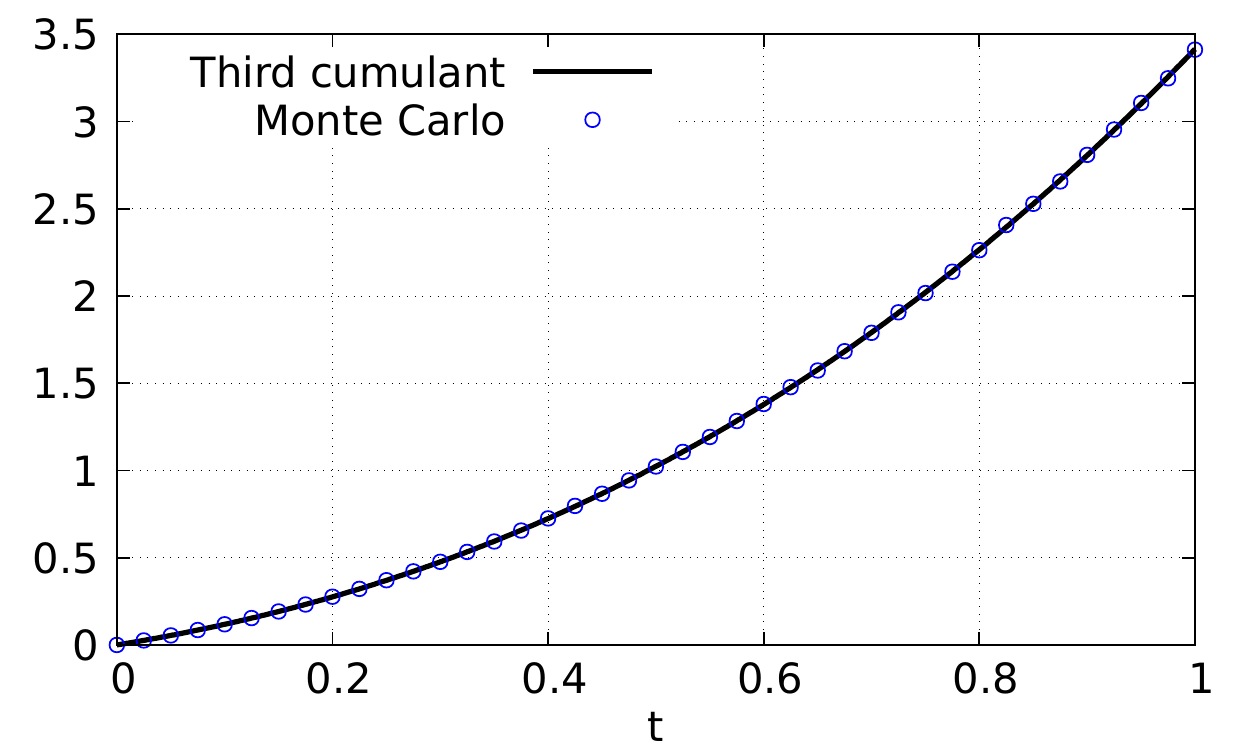}
    \caption{Third cumulant $\kappa^{(3)} (t)$.} 
   \end{subfigure}
  \begin{subfigure}[b]{0.49\textwidth}
    \includegraphics[width=1\linewidth, height=5cm]{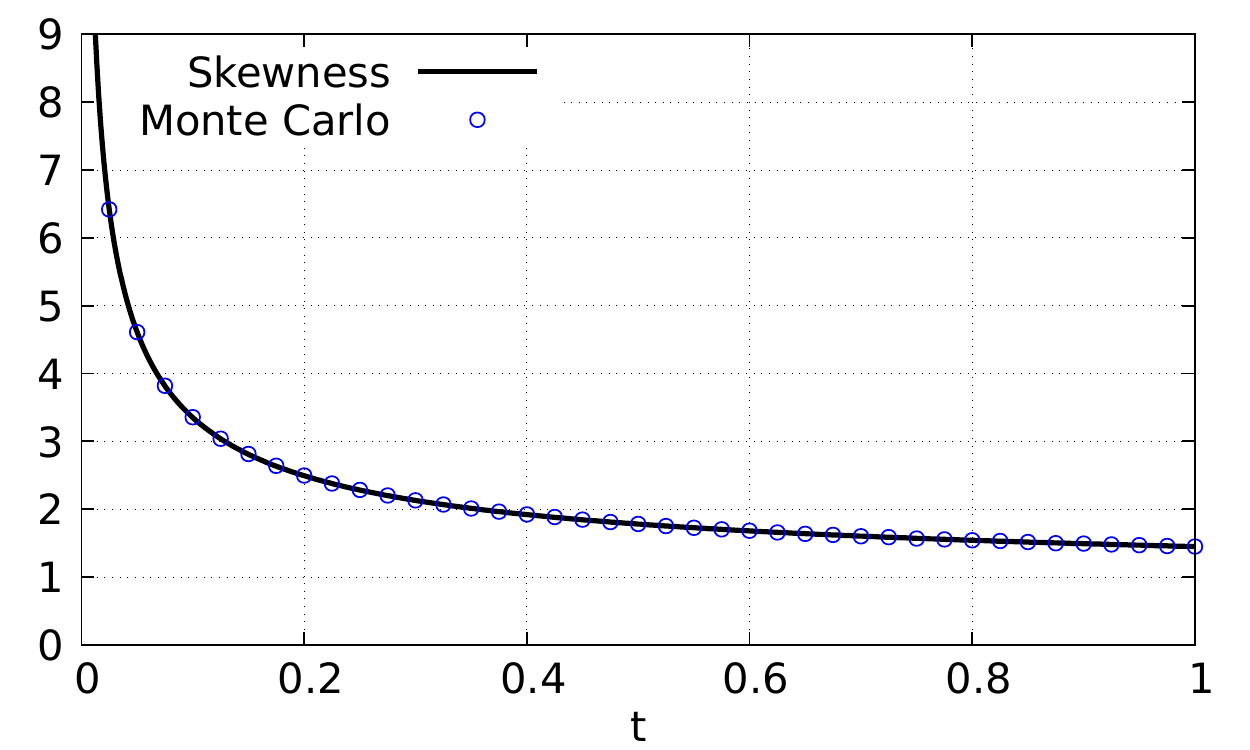} 
 \caption{Skewness $\kappa^{(3)} (t)/\big( \kappa^{(2)} (t)\big)^{3/2}$.} 
 \end{subfigure}
   \caption{Third cumulant and skewness with exponential kernel.} 
 \label{fig2} 
\end{figure}

\noindent 
Figure~\ref{fig2} shows the numerical evaluation of
$\kappa^{(3)}(t)$ and of the associated skewness
$\kappa^{(3)}(t)/ (\kappa^{(2)}(t))^{3/2}$. 

\subsubsection*{Fourth cumulant} 
 The recursive computation of $\kappa^{(4)} (t)$ can be similarly carried out
 from \eqref{k4} and \eqref{c1jk}-\eqref{c1jk2} using Mathematica, which yields 
\begin{align*} 
  & 
  \kappa^{(4)}(t) =
  -\frac{\nu }{12 (a-b)^8}
\Big(
180 a b^6
+570 a^2 b^5
+100 a^3 b^4
-15 a^4 b^3
+2 a^5 b^2
+
12b^5 (
6 a^3
+2 a^2 b
- b^3
-7 a b^2
) t
\\
&
+
4 \big( 
5 a^6 b
-45 a b^6
+3 a^7
-59a^5 b^2
-180 a^2 b^5
+75 a^3 b^4
+75 a^4 b^3
\\
&
+
6 a^2 b (
5 a b^4
-25 b^5
+41 a^2 b^3
-22 a^3 b^2
- a^4 b
+2 a^5 ) t 
+ 18 a^3 b^2 ( 
 10 a b^3
-5 b^4
-4 a^2 b^2
-2 a^3 b
+ a^4 ) t^2 
\\
&
+ 12 a^4b^3 (a-b)^3
 t^3 
\big) \re^{(a-b) t}
+
\big( 
150 a^2 b^5
-360 a^3 b^4
-564 a^4 b^3
+588 a^5 b^2
+18 a^6 b
-84 a^7
\\
&
+ 4 a^3 b ( 
90 b^4
-306 a b^3
+180 a^2 b^2
+108 a^3 b
-72 a^4 ) t 
+ 144 a^4 b^2 ( 
- 4 a b^2
+ b^3
+5 a^2 b
- 2 a^3 
)t^2 \big)
\re^{2 (a-b) t}
\\
&
+ \big( 
 276 a^4 b^3
-40 a^3 b^4
-132 a^6 b
-320 a^5 b^2
+144 a^7
+ 24 a^4 b (
13 a b^2
-23 a^2 b
-2 b^3
+12 a^3 )t 
\big) \re^{3 (a-b) t}
\\
&
+ a^4 \big( 
3 b^3
-34 a b^2
+94 a^2 b
-72 a^3 \big) \re^{4 (a-b) t}
\Big).
\end{align*} 

\noindent
Figure~\ref{fig2} shows the numerical evaluation of
$\kappa^{(4)}(t)$ and of the associated excess kurtosis 
$\kappa^{(4)}(t)/ (\kappa^{(2)}(t))^2$. 

\begin{figure}[H]
  \centering
 \begin{subfigure}[b]{0.49\textwidth}
    \includegraphics[width=1\linewidth, height=5cm]{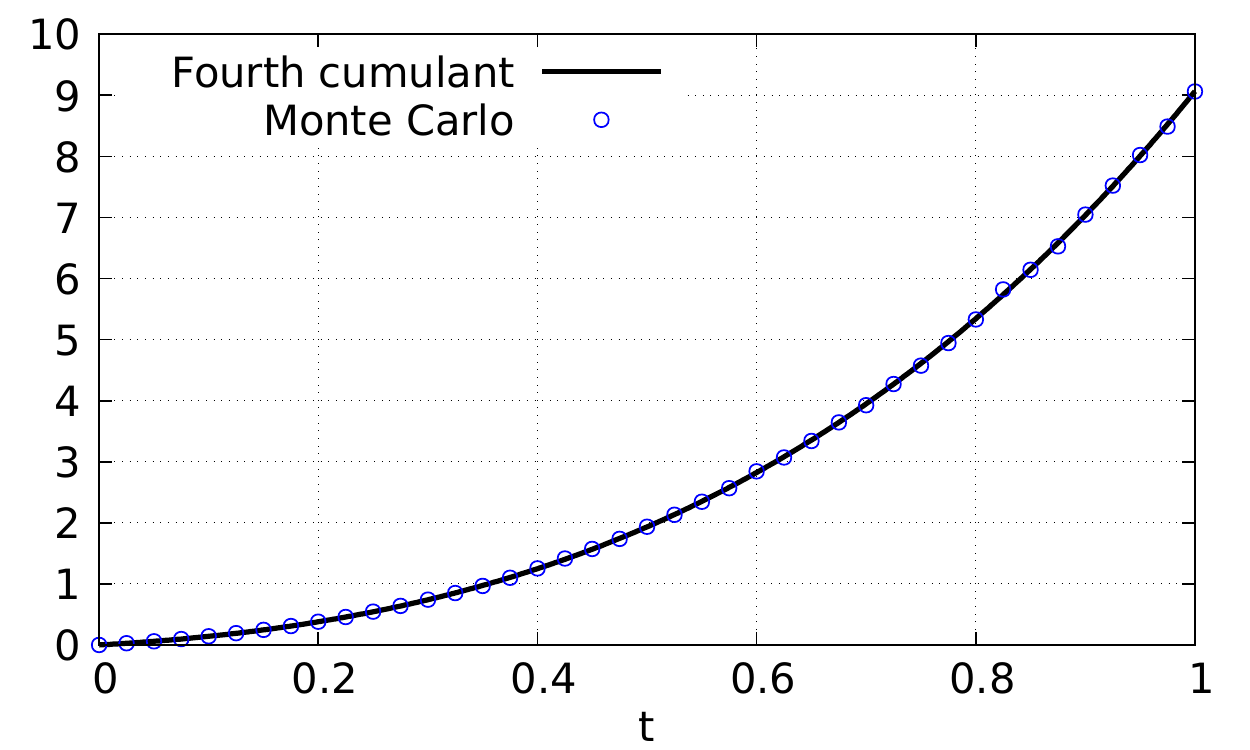} 
    \caption{Fourth cumulant $\kappa^{(4)} (t)$.} 
    \end{subfigure}
  \begin{subfigure}[b]{0.49\textwidth}
    \includegraphics[width=1\linewidth, height=5cm]{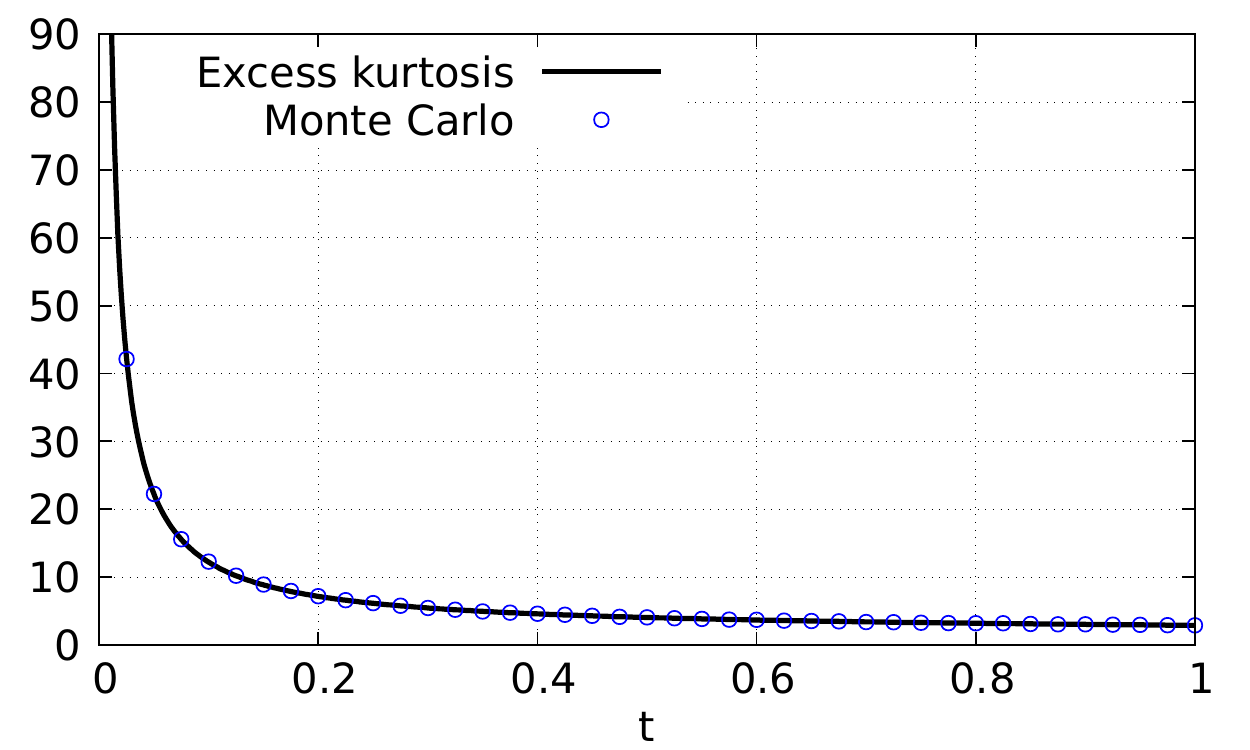} 
    \caption{Excess kurtosis $\kappa^{(4)} (t) / \big( \kappa^{(2)} (t) \big)$.} 
   \end{subfigure}
    \caption{Fourth cumulant and excess kurtosis with exponential kernel.} 
\end{figure}

\subsubsection*{Intensity cumulants} 
\noindent
 We have 
$$ 
 \E [ \lambda_t ] =
 a \E \left[ \int_0^t \re^{-b(t-s)} dN_s \right] 
 = a \re^{-bt} \kappa^{(1)}(e_{0,b}) = 
 a \re^{-bt} \int_0^t \kappa^{(1)}_z (e_{0,b}) \nu (dz),
 $$
 with
\begin{eqnarray*} 
 \kappa^{(1)}_z (e_{0,b}) 
 & = & e_{0,b}(z) 
 + \sum_{n=1}^\infty
 \int_0^t \cdots \int_0^t
 e_{0,b}(z+x_1+\cdots + x_n )
  \gamma ( dx_1) \cdots \gamma ( dx_n) 
  \\
  & = &
  \re^{bz} e_{0,a}(t-z)
    \\
  & = &
    \re^{bt} e_{0,a-b}(t-z)
      \\
  & = &
  \re^{at} e_{0,b-a}(z), 
\end{eqnarray*} 
 hence
$$ 
 \E [ \lambda_t ] =
 a \re^{-bt} \kappa^{(1)} (e_{0,b}) 
  = 
\nu a \re^{-bt} \int_0^t
\kappa^{(1)}_z (e_{0,b}) \nu (dz)
= 
\frac{\nu a}{b-a} ( 1-\re^{(a-b)t}), 
$$ 
see e.g. Theorem~3.6 in \cite{dassios-zhao2}.
 Next, we compute the joint moment
$\E [ \lambda_t N_t]$.
Using \eqref{c1jk}, we have
\begin{eqnarray*} 
\lefteqn{
 \kappa^{(2)}_z (e_{0,0},e_{0,b}) = 
 \sum_{m=1}^\infty
 \int_0^t \cdots \int_0^t
 \kappa_{z+x_1+\cdots + x_m}^{(1)}(e_{0,0})
 \kappa_{z+x_1+\cdots + x_m}^{(1)} (e_{0,b}) 
  \gamma ( dx_1) \cdots \gamma ( dx_m) 
  }
 \\
 & = & 
  \re^{at}
  \frac{b}{b-a}
  \sum_{m=1}^\infty
 \int_0^t \cdots \int_0^t
 e_{0,b-a}(z+x_1+\cdots + x_m)
 \gamma ( dx_1) \cdots \gamma ( dx_m) 
  \\
 & & - \frac{a}{b-a}
  \re^{(2a-b )t}
  \sum_{m=1}^\infty
 \int_0^t \cdots \int_0^t
 e_{0,2(b-a)} (z+x_1+\cdots + x_m) 
 \gamma ( dx_1) \cdots \gamma ( dx_m)
 \\
 & = & 
  \frac{ab}{b-a}
  \re^{at}
  z 
 + a^2 \re^{(2a-b )t}
 \frac{1 - \re^{(b-a)(t-z)} }{(b-a)^2}, \qquad z\in [0,t].
\end{eqnarray*} 
\noindent
 Hence we have 
\begin{eqnarray*} 
  \lefteqn{
    \E [ \lambda_t N_t]
    =
    a \E \left[ N_t \int_0^t \re^{-b(t-s)} dN_s \right] 
    = a \re^{-bt} \kappa^{(2)}(e_{0,0},e_{0,b}) 
+ a \re^{-bt} \kappa^{(1)} (e_{0,b})\kappa^{(1)} (e_{0,0})
  }
  \\
& = & 
a \re^{-bt} \int_0^t \kappa_z^{(2)}(e_{0,0},e_{0,b}) \nu ( dz)
  +
a \re^{-bt}  \int_0^t \kappa_z^{(1)}(e_{0,0}) \kappa_z^{(1)}(e_{0,b}) 
\nu ( dz)
+ a \re^{-bt} \kappa^{(1)} (e_{0,b})\kappa^{(1)} (e_{0,0})
\\
& = & 
 -\frac{\nu a }{(a - b)^3} 
  \Big( b^2 - \nu a - a b 
  +
  \nu b ( b -a ) t
  + a \big( a - \nu + a ( b - a ) t \big) \re^{2 (a-b) t}
  \\
 &  &  +
  \big( 2 \nu a - a^2 + a b - b^2 
  +
    ( \nu a b - \nu b^2 -a ( a-b)^2 ) t
 +
  ab ( a-b)^2 t^2 /2 \big) \re^{(a - b) t} 
    \Big), \quad t\geq 0. 
\end{eqnarray*}

\noindent 
 The following figures are plotted with $\nu=2$, $a=0.5$, $b=1$, and $10^6$
 Monte Carlo samples.
 
\begin{figure}[H]
  \centering
 \begin{subfigure}[b]{0.49\textwidth}
    \includegraphics[width=1\linewidth, height=5cm]{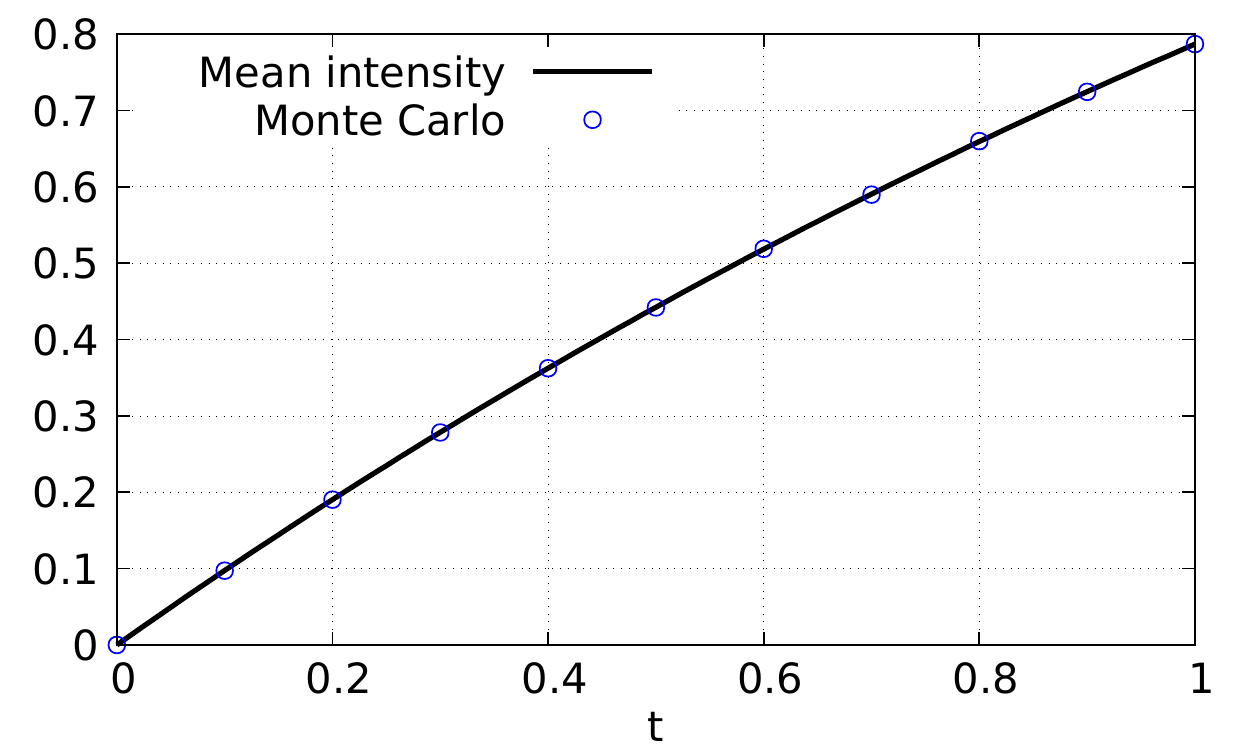} 
    \caption{Mean intensity $\E [ \lambda_t ]$.} 
   \end{subfigure}
  \begin{subfigure}[b]{0.49\textwidth}
    \includegraphics[width=1\linewidth, height=5cm]{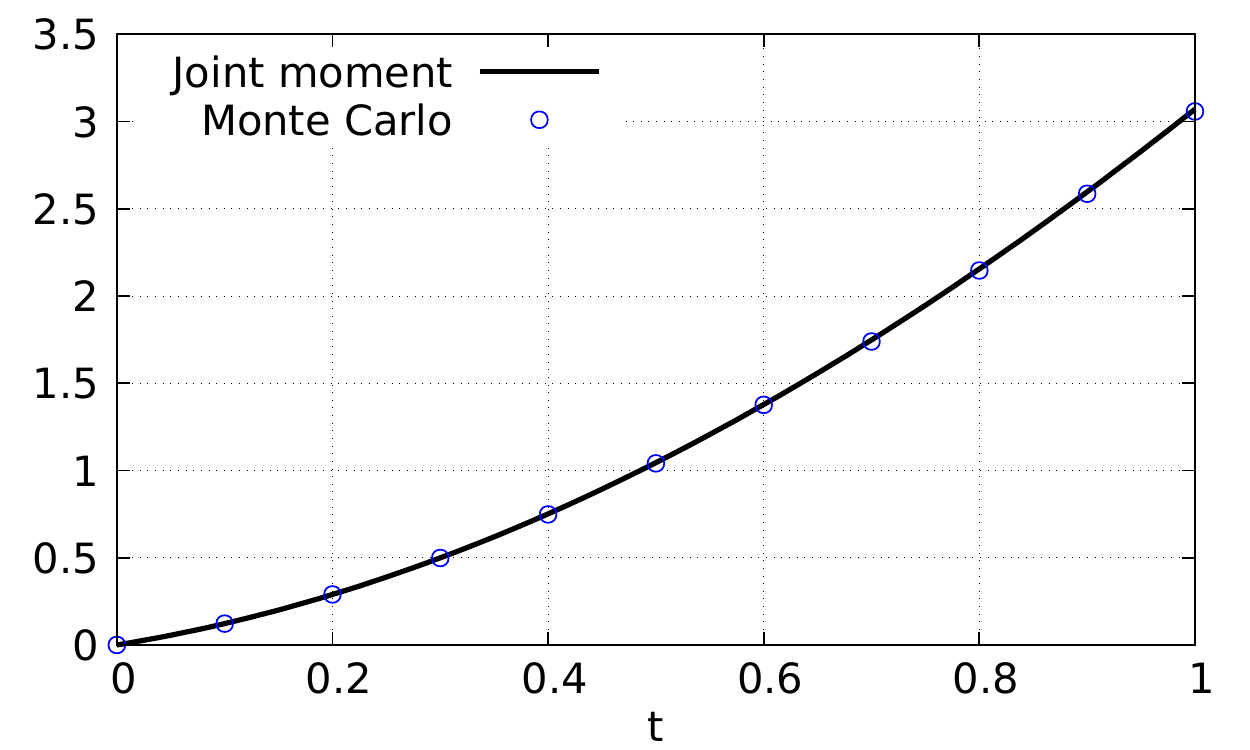} 
    \caption{Joint moment $\E [ \lambda_t N_t]$.} 
    \end{subfigure}
    \caption{Mean intensity and joint moment with exponential kernel.} 
\end{figure}
 
\footnotesize 

\def\cprime{$'$} \def\polhk#1{\setbox0=\hbox{#1}{\ooalign{\hidewidth
  \lower1.5ex\hbox{`}\hidewidth\crcr\unhbox0}}}
  \def\polhk#1{\setbox0=\hbox{#1}{\ooalign{\hidewidth
  \lower1.5ex\hbox{`}\hidewidth\crcr\unhbox0}}} \def\cprime{$'$}

\end{document}